\documentclass[12pt,reqno]{amsart} 
\usepackage{amssymb}
\begin{document}
\theoremstyle{plain}
\newtheorem{thm}{Theorem}[section]
\newtheorem{theorem}[thm]{Theorem}
\newtheorem{lemma}[thm]{Lemma}
\newtheorem{corollary}[thm]{Corollary}
\newtheorem{proposition}[thm]{Proposition}
\newtheorem{example}[thm]{Example}
\theoremstyle{definition}
\newtheorem{notation}[thm]{Notation}
\newtheorem{claim}[thm]{Claim}
\newtheorem{remark}[thm]{Remark}
\newtheorem{remarks}[thm]{Remarks}
\newtheorem{conjecture}[thm]{Conjecture}
\newtheorem{definition}[thm]{Definition}
\newtheorem{problem}[thm]{Problem}
\newcommand{\fz}{\frak{z}}
\newcommand{\zar}{{\rm zar}}
\newcommand{\an}{{\rm an}} 
\newcommand{\red}{{\rm red}}
\newcommand{\codim}{{\rm codim}}
\newcommand{\rank}{{\rm rank}}
\newcommand{\Pic}{{\rm Pic}}
\newcommand{\Div}{{\rm Div}}
\newcommand{\Hom}{{\rm Hom}}
\newcommand{\im}{{\rm im}}
\newcommand{\Spec}{{\rm Spec}}
\newcommand{\sing}{{\rm sing}}
\newcommand{\reg}{{\rm reg}}
\newcommand{\Char}{{\rm char}}
\newcommand{\Tr}{{\rm Tr}}
\newcommand{\Gal}{{\rm Gal}}
\newcommand{\Min}{{\rm Min \ }}
\newcommand{\Max}{{\rm Max \ }}
\newcommand{\soplus}[1]{\stackrel{#1}{\oplus}}
\newcommand{\dlog}{{\rm dlog}\,}    
\newcommand{\limdir}[1]{{\displaystyle{\mathop{\rm
lim}_{\buildrel\longrightarrow\over{#1}}}}\,} 
\newcommand{\liminv}[1]{{\displaystyle{\mathop{\rm
lim}_{\buildrel\longleftarrow\over{#1}}}}\,} 
\newcommand{\boxtensor}{{\Box\kern-9.03pt\raise1.42pt\hbox{$\times$}}}
\newcommand{\sext}{\mbox{${\mathcal E}xt\,$}}
\newcommand{\shom}{\mbox{${\mathcal H}om\,$}}
\newcommand{\coker}{{\rm coker}\,} 
\renewcommand{\iff}{\mbox{ $\Longleftrightarrow$ }}
\newcommand{\onto}{\mbox{$\,\>>>\hspace{-.5cm}\to\hspace{.15cm}$}}
\catcode`\@=11
\def\opn#1#2{\def#1{\mathop{\kern0pt\fam0#2}\nolimits}} 
\def\bold#1{{\bf #1}}%
\def\underrightarrow{\mathpalette\underrightarrow@}
\def\underrightarrow@#1#2{\vtop{\ialign{$##$\cr
 \hfil#1#2\hfil\cr\noalign{\nointerlineskip}%
 #1{-}\mkern-6mu\cleaders\hbox{$#1\mkern-2mu{-}\mkern-2mu$}\hfill
 \mkern-6mu{\to}\cr}}}
\let\underarrow\underrightarrow
\def\underleftarrow{\mathpalette\underleftarrow@}
\def\underleftarrow@#1#2{\vtop{\ialign{$##$\cr
 \hfil#1#2\hfil\cr\noalign{\nointerlineskip}#1{\leftarrow}\mkern-6mu
 \cleaders\hbox{$#1\mkern-2mu{-}\mkern-2mu$}\hfill
 \mkern-6mu{-}\cr}}}
\let\amp@rs@nd@\relax
\newdimen\ex@
\ex@.2326ex
\newdimen\bigaw@
\newdimen\minaw@
\minaw@16.08739\ex@
\newdimen\minCDaw@
\minCDaw@2.5pc
\newif\ifCD@
\def\minCDarrowwidth#1{\minCDaw@#1}
\newenvironment{CD}{\@CD}{\@endCD}
\def\@CD{\def\A##1A##2A{\llap{$\vcenter{\hbox
 {$\scriptstyle##1$}}$}\Big\uparrow\rlap{$\vcenter{\hbox{%
$\scriptstyle##2$}}$}&&}%
\def\V##1V##2V{\llap{$\vcenter{\hbox
 {$\scriptstyle##1$}}$}\Big\downarrow\rlap{$\vcenter{\hbox{%
$\scriptstyle##2$}}$}&&}%
\def\={&\hskip.5em\mathrel
 {\vbox{\hrule width\minCDaw@\vskip3\ex@\hrule width
 \minCDaw@}}\hskip.5em&}%
\def\verteq{\Big\Vert&&}%
\def\noarr{&&}%
\def\vspace##1{\noalign{\vskip##1\relax}}\relax\let\amp@rs@nd@&\iffalse}\fi
 \CD@true\vcenter\bgroup\relax\let\\=\cr\iffalse}\fi\tabskip\z@skip\baselineskip20\ex@
 \lineskip3\ex@\lineskiplimit3\ex@\halign\bgroup
 &\hfill$\m@th##$\hfill\cr}
\def\@endCD{\cr\egroup\egroup}
\def\>#1>#2>{\amp@rs@nd@\setbox\z@\hbox{$\scriptstyle
 \;{#1}\;\;$}\setbox\@ne\hbox{$\scriptstyle\;{#2}\;\;$}\setbox\tw@
 \hbox{$#2$}\ifCD@
 \global\bigaw@\minCDaw@\else\global\bigaw@\minaw@\fi
 \ifdim\wd\z@>\bigaw@\global\bigaw@\wd\z@\fi
 \ifdim\wd\@ne>\bigaw@\global\bigaw@\wd\@ne\fi
 \ifCD@\hskip.5em\fi
 \ifdim\wd\tw@>\z@
 \mathrel{\mathop{\hbox to\bigaw@{\rightarrowfill}}\limits^{#1}_{#2}}\else
 \mathrel{\mathop{\hbox to\bigaw@{\rightarrowfill}}\limits^{#1}}\fi
 \ifCD@\hskip.5em\fi\amp@rs@nd@}
\def\<#1<#2<{\amp@rs@nd@\setbox\z@\hbox{$\scriptstyle
 \;\;{#1}\;$}\setbox\@ne\hbox{$\scriptstyle\;\;{#2}\;$}\setbox\tw@
 \hbox{$#2$}\ifCD@
 \global\bigaw@\minCDaw@\else\global\bigaw@\minaw@\fi
 \ifdim\wd\z@>\bigaw@\global\bigaw@\wd\z@\fi
 \ifdim\wd\@ne>\bigaw@\global\bigaw@\wd\@ne\fi
 \ifCD@\hskip.5em\fi
 \ifdim\wd\tw@>\z@
 \mathrel{\mathop{\hbox to\bigaw@{\leftarrowfill}}\limits^{#1}_{#2}}\else
 \mathrel{\mathop{\hbox to\bigaw@{\leftarrowfill}}\limits^{#1}}\fi
 \ifCD@\hskip.5em\fi\amp@rs@nd@}
\newenvironment{CDS}{\@CDS}{\@endCDS}
\def\@CDS{\def\A##1A##2A{\llap{$\vcenter{\hbox
 {$\scriptstyle##1$}}$}\Big\uparrow\rlap{$\vcenter{\hbox{%
$\scriptstyle##2$}}$}&}%
\def\V##1V##2V{\llap{$\vcenter{\hbox
 {$\scriptstyle##1$}}$}\Big\downarrow\rlap{$\vcenter{\hbox{%
$\scriptstyle##2$}}$}&}%
\def\={&\hskip.5em\mathrel
 {\vbox{\hrule width\minCDaw@\vskip3\ex@\hrule width
 \minCDaw@}}\hskip.5em&}
\def\verteq{\Big\Vert&}
\def\novarr{&}
\def\noharr{&&}
\def\SE##1E##2E{\slantedarrow(0,18)(4,-3){##1}{##2}&}
\def\SW##1W##2W{\slantedarrow(24,18)(-4,-3){##1}{##2}&}
\def\NE##1E##2E{\slantedarrow(0,0)(4,3){##1}{##2}&}
\def\NW##1W##2W{\slantedarrow(24,0)(-4,3){##1}{##2}&}
\def\slantedarrow(##1)(##2)##3##4{%
\thinlines\unitlength1pt\lower 6.5pt\hbox{\begin{picture}(24,18)%
\put(##1){\vector(##2){24}}%
\put(0,8){$\scriptstyle##3$}%
\put(20,8){$\scriptstyle##4$}%
\end{picture}}}
\def\vspace##1{\noalign{\vskip##1\relax}}\relax\let\amp@rs@nd@&\iffalse}\fi
 \CD@true\vcenter\bgroup\relax\let\\=\cr\iffalse}\fi\tabskip\z@skip\baselineskip20\ex@
 \lineskip3\ex@\lineskiplimit3\ex@\halign\bgroup
 &\hfill$\m@th##$\hfill\cr}
\def\@endCDS{\cr\egroup\egroup}
\newdimen\TriCDarrw@
\newif\ifTriV@
\newenvironment{TriCDV}{\@TriCDV}{\@endTriCD}
\newenvironment{TriCDA}{\@TriCDA}{\@endTriCD}
\def\@TriCDV{\TriV@true\def\TriCDpos@{6}\@TriCD}
\def\@TriCDA{\TriV@false\def\TriCDpos@{10}\@TriCD}
\def\@TriCD#1#2#3#4#5#6{%
\setbox0\hbox{$\ifTriV@#6\else#1\fi$}
\TriCDarrw@=\wd0 \advance\TriCDarrw@ 24pt
\advance\TriCDarrw@ -1em
\def\SE##1E##2E{\slantedarrow(0,18)(2,-3){##1}{##2}&}
\def\SW##1W##2W{\slantedarrow(12,18)(-2,-3){##1}{##2}&}
\def\NE##1E##2E{\slantedarrow(0,0)(2,3){##1}{##2}&}
\def\NW##1W##2W{\slantedarrow(12,0)(-2,3){##1}{##2}&}
\def\slantedarrow(##1)(##2)##3##4{\thinlines\unitlength1pt
\lower 6.5pt\hbox{\begin{picture}(12,18)%
\put(##1){\vector(##2){12}}%
\put(-4,\TriCDpos@){$\scriptstyle##3$}%
\put(12,\TriCDpos@){$\scriptstyle##4$}%
\end{picture}}}
\def\={\mathrel {\vbox{\hrule
   width\TriCDarrw@\vskip3\ex@\hrule width
   \TriCDarrw@}}}
\def\>##1>>{\setbox\z@\hbox{$\scriptstyle
 \;{##1}\;\;$}\global\bigaw@\TriCDarrw@
 \ifdim\wd\z@>\bigaw@\global\bigaw@\wd\z@\fi
 \hskip.5em
 \mathrel{\mathop{\hbox to \TriCDarrw@
{\rightarrowfill}}\limits^{##1}}
 \hskip.5em}
\def\<##1<<{\setbox\z@\hbox{$\scriptstyle
 \;{##1}\;\;$}\global\bigaw@\TriCDarrw@
 \ifdim\wd\z@>\bigaw@\global\bigaw@\wd\z@\fi
 \mathrel{\mathop{\hbox to\bigaw@{\leftarrowfill}}\limits^{##1}}
 }
 \CD@true\vcenter\bgroup\relax\let\\=\cr\iffalse}\fi
 \tabskip\z@skip\baselineskip20\ex@
 \lineskip3\ex@\lineskiplimit3\ex@
 \ifTriV@
 \halign\bgroup
 &\hfill$\m@th##$\hfill\cr
#1&\multispan3\hfill$#2$\hfill&#3\\
&#4&#5\\
&&#6\cr\egroup%
\else
 \halign\bgroup
 &\hfill$\m@th##$\hfill\cr
&&#1\\%
&#2&#3\\
#4&\multispan3\hfill$#5$\hfill&#6\cr\egroup
\fi}
\def\@endTriCD{\egroup}
\newcommand{\sA}{{\mathcal A}}
\newcommand{\sB}{{\mathcal B}}
\newcommand{\sC}{{\mathcal C}}
\newcommand{\sD}{{\mathcal D}}
\newcommand{\sE}{{\mathcal E}}
\newcommand{\sF}{{\mathcal F}}
\newcommand{\sG}{{\mathcal G}}
\newcommand{\sH}{{\mathcal H}}
\newcommand{\sI}{{\mathcal I}}
\newcommand{\sJ}{{\mathcal J}}
\newcommand{\sK}{{\mathcal K}}
\newcommand{\sL}{{\mathcal L}}
\newcommand{\sM}{{\mathcal M}}
\newcommand{\sN}{{\mathcal N}}
\newcommand{\sO}{{\mathcal O}}
\newcommand{\sP}{{\mathcal P}}
\newcommand{\sQ}{{\mathcal Q}}
\newcommand{\sR}{{\mathcal R}}
\newcommand{\sS}{{\mathcal S}}
\newcommand{\sT}{{\mathcal T}}
\newcommand{\sU}{{\mathcal U}}
\newcommand{\sV}{{\mathcal V}}
\newcommand{\sW}{{\mathcal W}}
\newcommand{\sX}{{\mathcal X}}
\newcommand{\sY}{{\mathcal Y}}
\newcommand{\sZ}{{\mathcal Z}}
\newcommand{\A}{{\mathbb A}}
\newcommand{\B}{{\mathbb B}}
\newcommand{\C}{{\mathbb C}}
\newcommand{\D}{{\mathbb D}}
\newcommand{\E}{{\mathbb E}}
\newcommand{\F}{{\mathbb F}}
\newcommand{\G}{{\mathbb G}}
\newcommand{\HH}{{\mathbb H}}
\newcommand{\I}{{\mathbb I}}
\newcommand{\J}{{\mathbb J}}
\newcommand{\M}{{\mathbb M}}
\newcommand{\N}{{\mathbb N}}
\renewcommand{\P}{{\mathbb P}}
\newcommand{\Q}{{\mathbb Q}}
\newcommand{\R}{{\mathbb R}}
\newcommand{\T}{{\mathbb T}}
\newcommand{\U}{{\mathbb U}}
\newcommand{\V}{{\mathbb V}}
\newcommand{\W}{{\mathbb W}}
\newcommand{\X}{{\mathbb X}}
\newcommand{\Y}{{\mathbb Y}}
\newcommand{\Z}{{\mathbb Z}}
\title[Transitive and Self-dual Codes]{ Transitive and self-dual codes attaining the 
Tsfasman-Vladut-Zink bound}   
\author{Henning Stichtenoth}  
\address{Department of Mathematics, University of Duisburg-Essen, Campus Essen, 
45117 Essen, Germany}
\email{stichtenoth@uni-essen.de}
\address{Sabanci University, MDBF, 34956 Orhanli, Tuzla, Istanbul, Turkey}
\email{henning@sabanciuniv.edu}
\begin{abstract}
A major problem in coding theory is the question if the class of cyclic codes is asymptotically good.
 In this paper we introduce - as a generalization of cyclic codes - the notion of {\it transitive} 
codes (see Definition 1.4 in Section 1),
 and we show that the class of transitive codes is asymptotically good.
 Even more, transitive codes attain the Tsfasman-Vladut-Zink bound over $\F_q$, for all
 squares $q = l^2$. We also show that self-orthogonal and self-dual codes attain the 
Tsfasman-Vladut-Zink bound, thus improving previous results about self-dual codes
 attaining the Gilbert-Varshamov bound. The main tool is a new asymptotically optimal 
tower $E_0 \subseteq E_1 \subseteq E_2 \subseteq \ldots$ of function fields over
 $\F_q$ (with $q = \ell ^2)$, where all extensions $E_n /E_0$ are Galois. 
\end{abstract}

\maketitle


\section{Introduction and Main Results}

Let $\F_q$ be the finite field of cardinality $q$. In this paper we consider primarily {\it linear} 
$[n,k,d]$-codes  $C$ over $\F_q$; i.e., the parameters $n = n (C),     \ k = k (C)$ and $d = d (C)$ 
are the {\it length}, the {\it dimension} and the {\it minimum distance} of the code. The ratios 
$R= R(C) = k(C)/n(C)$ and $\delta = \delta (C) = d (C)/n(C)$ denote the 
{\it information rate} and the {\it relative minimum distance}, resp., of the code. \\

A crucial role in the asymptotic theory of codes plays the set $U_q \subseteq [0,1] \times [0,1]$ 
which is defined as follows: a point $(\delta , R) \in \R^2$ with $0 \leq \delta \leq 1$ and
 $0 \leq R \leq 1$ belongs to $U_q$ if and only if there exists a sequence $(C_i)_{i \geq 0}$
 of codes over $\F_q$ such that 
$$
n (C_i) \to \infty, \ \delta (C_i) \to \delta \ \mbox{and} \ R(C_i) \to R, \ \mbox{as} \ i \to \infty. 
$$
One then defines the function $\alpha_q : [0,1] \to [0,1]$ by 
$$
\alpha_q (\delta) = \sup \{ R ; (\delta ,R) \in U_q \} , \ \mbox{\rm for} \  \delta \in [0,1].  
$$
The following facts are well-known (and easy to prove), see [M], [T-V]: 

\begin{proposition}
{\rm i)} A point $(\delta, R) \in [0,1] \times [0,1]$ belongs to the set $U_q$ if and only 
if  \ $0 \leq R \leq \alpha_q (\delta)$. \\
{\rm ii)} The function $\alpha_q$ is continuous and non-increasing.\\
{\rm iii)} $\alpha_q (0) =1$, and $\alpha_q (\delta) =0$ for $1 - q^{-1} \leq \delta \leq 1$.  
\end{proposition}

Many {\it upper bounds} for $\alpha_q (\delta)$ are known, see [vL], [T-V]. More 
interesting, however, are {\it lower bounds} for $\alpha_q (\delta)$, since any 
non-trivial lower bound assures the existence of arbritrarily long linear codes 
with good error correction parameters. The classical lower bound for
 $\alpha_q (\delta)$ is the {\it asymptotic Gilbert-Varshamov bound}, which says:

\begin{proposition} {\rm (see [vL])}. For all $\delta \in (0,1 - q^{-1} )$ one has 
$$
\alpha_q (\delta) \geq 1 - \delta \log_q (q-1) + \delta \log_q (\delta) + (1- \delta) \log_q (1 - \delta).
$$
\end{proposition}

For sufficiently large non-prime $q$ and for certain ranges of the variable $\delta$, 
the Gilbert-Varshamov bound is improved by the {\it Tsfasman-Vladut-Zink bound}
 as follows: 

\begin{proposition} {\rm (see [T-V-Z], [N-X])}. Let 
$$
A (q) = \lim\sup_{g \to \infty} N_q (g)/g ,
$$
where $N_q (g)$ denotes the maximum number of rational places that a function fields $F$ 
over $\F_q$ of genus $g$ can have. Then
$$
\alpha_q (\delta) \geq 1 - \delta - A (q)^{-1} \ \ \mbox{for} \ 0 \leq \delta \leq 1. 
$$
\end{proposition}

It is well known that $A (q) \leq q^{1/2} -1$ (this is the {\it Drinfeld-Vladut bound}), and
 $A (q) = q^{1/2} -1 $ if $q$ is a square, see [I], [T-V-Z], [G-S1]. It then follows easily
 that the Tsfasman-Vladut-Zink bound in Proposition 1.3 improves the Gilbert-Varshamov bound
 for all squares $q \geq 49$. 
For {\it non-linear codes} over $\F_q$, the Tsfasman-Vladut-Zink bound was further 
improved recently, see [X], [N-O1], [N-O2],  [S-X], [E1], [E2].\\

 In order to prove the Gilbert-Varshamov
 and the Tsfasman-Vladut-Zink bound one constructs families of long codes with sufficiently 
good parameters. However, the proofs provide linear codes without any particular
structure. For instance, one of the most challenging problems in coding theory is
 still open (see [P-H], [M-W]): Do there exist sequences $(C_i)_{i \geq 0}$
 of {\it cyclic codes} $C_i$ over $\F_q$ with 
$$
n (C_i) \to \infty, \ \  \lim_{i \to \infty} R(C_i) > 0 \ \ \mbox{and} \ \lim_{i \to \infty} \delta (C_i) > 0 \ ?
$$
Cyclic codes can be understood as a special case of what we call in this paper
 {\it transitive codes}. Recall that a subgroup $U$ of the symmetric group $S_n$
 is called {\it transitive} if for any pair $(i,j)$ with $i,j \in \{ 1, \ldots ,n\}$ there is a
 permutation $\pi \in U$ such that $\pi (i) = j$. A permutation $\pi \in S_n$ is called 
an {\it automorphism} of the code $C \subseteq \F^{n}_{q}$ if 
$$
(c_1, \ldots , c_n) \in C \Rightarrow (c_{\pi (1)}, \ldots , c_{\pi (n)}) \in C
$$
holds for all codewords $(c_1, \ldots , c_n) \in C$. The {\it automorphism group}
 ${\rm Aut} (C) \subseteq S_n$ is the group of all automorphisms of the code $C$. 

\begin{definition}
A code $C$ over $\F_q$ of length $n$ is said to be {\it transitive} if its automorphism 
group ${\rm Aut} (C)$ is a transitive subgroup of $S_n$. 
\end{definition}

It is obvious that any cyclic code is transitive. We can now state our first result. 

\begin{theorem}
Let $q = \ell ^2$ be a square. Then the class of transitive codes meets the
 Tsfasman-Vladut-Zink bound. More precisely, let $R, \delta \geq 0$ be real numbers with
 $R = 1 - \delta - 1/ (\ell -1)$. Then there exists a sequence $(C_j)_{j \geq 0}$ 
of linear codes $C_j$ over $\F_q$ with parameters $[n_j, k_j, d_j]$ with the 
following properties: 
\begin{enumerate}
\item[a)] All $C_j$ are transitive codes.
\item[b)] $n_j \to \infty$ as $j \to \infty$. 
\item[c)] $lim_{j \to \infty} k_j / n_j \geq R$ and $\lim_{j \to \infty} d_j /n_j \geq \delta$. 
\end{enumerate}
\end{theorem}

Other important classes of codes are the {\it self-orthogonal codes} and the {\it self-dual 
codes}. Recall that a linear code $C$ is called self-orthogonal if $C$ is contained in its
 dual code $C^{\perp}$, and $C$ is called self-dual if $C = C^{\perp}$. It is clear that the
 information rate of a self-orthogonal codes satisfies $R(C) \leq 1/2$; the information 
rate of self-dual codes is $R(C) = 1/2$. It is well-known that self-dual codes reach the 
Gilbert-Varshamov bound, see [MW-S]. 
In this paper we shall prove: 

\begin{theorem}
Let $q = \ell ^2$ be a square. Then the class of self-orthogonal codes and the class of
 self-dual codes meet the Tsfasman-Vladut-Zink bound. More precisely we have: 
\begin{enumerate}
\item[i)] Let $0 \leq R \leq 1/2$ and $\delta \geq 0$ with $R = 1 - \delta - 1 /(\ell -1)$. Then
 there is a sequence $(C_j)_{j \geq 0}$ of linear codes $C_j$ over $\F_q$ with 
parameters $[n_j , k_j , d_j]$ such that: 
\begin{enumerate}
\item[a)] All $C_j$ are self-orthogonal codes. 
\item[b)] $n_j \to \infty$ as $j \to \infty$. 
\item[c)] $\lim_{j \to \infty} k_j / n_j \geq R$ and $\lim_{j \to \infty} d_j / n_j \geq \delta$. 
\end{enumerate}
\item[ii)] There is a sequence $(C_j)_{j \geq 0}$ of self-dual codes $C_j$ over $\F_q$ with
 parameters $[n_j , n_j/2, d_j]$ such that $n_j \to \infty$ and 
$$
\lim_{j \to \infty} d_j /n_j \geq 1/2 - 1 /(\ell -1). 
$$
\end{enumerate}
\end{theorem}

Note that the bounds given in Theorem 1.5 and Theorem 1.6 are better than the Gilbert-Varshamov
 bound, for all squares $q = \ell ^2 \geq 49$. \\

The main tool to prove Theorem 1.5 and Theorem 1.6 is a
 new asymptotically good tower of function fields over $\F_q$ which has particularly nice 
properties, see Theorem 1.7 below. Using that tower, we shall construct sequences of codes 
over $\F_q$ with the desired properties, analogously to the proof of Proposition 1.3 by 
Tsfasman-Vladut-Zink. \\

Before stating Theorem 1.7, we recall some notations from the theory of algebraic function 
fields, cf. [S1]. 
\begin{enumerate}
\item[-] For a function field $F/\F_q$ we denote by $g(F)$ the genus and by $N(F)$ the 
number of rational places of $F$. For an element $u \in F \setminus \{ 0 \}$, we denote
 by $(u)^F$, $(u)^{F}_{0}$ and $(u)^{F}_{\infty}$ the principal divisor, the zero divisor
 and the pole divisor, resp., of the element $u$. In particular we have $(u)^F = (u)^{F}_{0} - (u)^{F}_{\infty}$.
 The divisor of a differential $\mu \neq 0$ of $F/\F_q$ is denoted by $(\mu)^F$. 
\item[-] Let $\F_q (x)$ be a rational function field; then we denote, for $\alpha \in \F_q$,
 by $(x = \alpha)$ the zero of the function $(x - \alpha)$ and by $(x = \infty)$ the pole of 
the function $x$ in $\F_q (x)$. 
\item[-] Let $E/F$ be an extension of function fields over $\F_q$. Let $P$ be a place of 
$F$ and let $Q$ be a place of $E$ lying above $P$. Then $e (Q|P)$ and $d(Q|P)$ denote
 the ramification index and the different exponent, resp., of $Q|P$. The different of $E/F$ 
(which is a divisor of the function field $E$) is denoted by ${\rm Diff} (E/F)$. 
\end{enumerate}

\begin{theorem}
Let $q = \ell ^2$ be a square. Then there exists an infinite tower $\sE = 
(E_0 \subseteq E_1 \subseteq E_2 \subseteq \ldots)$ of function fields 
$E_i/\F_q$ with the following properties: 
\begin{enumerate}
\item[a)]$\F_q$ is the full constant field of $E_i$, for all $i \geq 0$. 
\item[b)] $E_0 = \F_q (z)$ is the rational function field. 
\item[c)] There exists an element $w \in E_1$ such that $w^{\ell -1} =z$. So
 we have $E_0 = \F_q (z) \subseteq \F_q (w) \subseteq E_1$, and the 
extension $\F_q (w) / E_0$ is cyclic of degree $(\ell -1)$. 
\item[d)] All extensions $E_n / E_0$ are Galois, and the degree of $E_n / E_0$ is 
$$
[E_n : E_0] = (\ell -1) \cdot \ell ^n \cdot p^{t(n)}, 
$$
where $p = {\rm char} (\F_q)$ is the characteristic of $\F_q$ and $t(n)$ 
is a non-negative integer. 
\item[e)] The place $(z=1)$ of $E_0$ splits completely in all extensions
 $E_n/E_0$; i.e., there are $[E_n: E_0]$ distinct places of $E_n$ above
 the place $(z=1)$, and all of them are rational places of $E_n$. In particular
 we have that the number of rational places satisfies
$N(E_n) \geq [E_n : E_0] = (\ell -1) \cdot \ell ^n \cdot p^{t(n)}$. 
\item[f)] The principal divisor of the function $w$ (as in item {\rm c)}) in the field 
$E_n$ has the form 
$$
(w)^{E_n} = e^{(n)}_{0} \cdot A^{(n)} - e^{(n)}_{\infty} \cdot B^{(n)}, 
$$
where $A^{(n)} > 0$ and $B^{(n)} > 0$ are positive divisors of the function field $E_n$. The
 ramification index $e^{(n)}_{0}$ of the place $(w =0)$ in $E_n / \F_q  (w)$ has the form 
$$
e^{(n)}_{0} = \ell ^{n-1} \cdot p^{r(n)} \ \mbox{with} \ r(n) \geq 0, 
$$
and the ramification index $e^{(n)}_{\infty}$ of the place $(w = \infty)$ in the extension 
$E_n /\F_q (w)$ has the form 
$$
e^{(n)}_{\infty} = \ell ^n \cdot p^{s(n)} \ \mbox{with} \ s (n) \geq 0.
$$
\item[g)] The different of the extension $E_n /\F_q (w)$ is given by 
$$
{\rm Diff} (E_n /\F_q (w)) = 2 (e^{(n)}_{0} -1)A^{(n)} + 2 (e^{(n)}_{\infty} -1) \cdot B^{(n)} , 
$$
with $e^{(n)}_{0}, e^{(n)}_{\infty}, A^{(n)}$ and $B^{(n)}$ as in item {\rm f)}. 
\item[h)] The genus $g (E_n)$ satisfies
$$
g (E_n)  = [E_n : \F_q (w)] +1 - (\deg A^{(n)} + \deg B^{(n)}) 
 \leq [E_n : \F_q (w)], 
 $$
with $A^{(n)}$ and $B^{(n)}$ as in item {\rm f)}. 
\item[j)] The tower $\sE$ attains the Drinfeld-Vladut bound; i.e., 
$$
\lim_{n \to \infty} N(E_n) / g(E_n) = q^{1/2} -1. 
$$
\end{enumerate}
\end{theorem}

This paper is organized as follows: In Section 2 we prove Theorem 1.7 which is the basis
 for our code constructions. In Section 3 we deal with transitive codes and give the proof 
of Theorem 1.5. We also explain briefly that the method of proof of Theorem 1.5 yields an
 improvement of the Tsfasman-Vladut-Zink bound for {\it transitive non-linear} codes.
 Finally, in Section 4 we discuss self-orthogonal and self-dual codes and we prove Theorem 1.6. 

\section{An Asymptotically Optimal Galois Tower of Function Fields}

For basic notations and facts in the theory of algebraic function fields we refer to [S1] and [N-X].
 We will in particular use the notations introduced in Section 1 after Theorem 1.6.\\

A {\it tower of function fields} over $\F_q$ is an infinite sequence 
$\sF = ( F_0, F_1, F_2, \ldots)$ of function fields $F_i$ over $\F_q$
 with the following properties:
\begin{enumerate}
\item[i)] $F_0 \subseteq F_1 \subseteq F_2 \subseteq \ldots $, and all
 extensions $F_{i+1}/F_i$ are separable of degree $[F_{i+1} : F_i] > 1$. 
\item[ii)] $\F_q$ is the full constant field of $F_i$, for all $i \geq 0$. 
\item[iii)] The genus $g(F_i)$ tends to infinity as $i \to \infty$. 
\end{enumerate}
Recall that $N(F_i)$ denotes the number of rational places of $F_i$ over 
$\F_q$. It is well-known that the {\it limit of the tower} $\sF$, 
$$
\lambda (\sF) : = \lim_{i \to \infty} N(F_i)/g(F_i)
$$
does exist (see [G-S2]). As follows from the Drinfeld-Vladut bound (see Section 1), one has that 
$$
0 \leq \lambda (\sF) \leq A (q) \leq q^{1/2} -1. 
$$
The tower $\sF$ is said to be {\it asymptotically optimal} if $\lambda (\sF) = A(q)$. For 
$q = \ell ^2$ a square number we have that $A(q) = \ell -1$, see Section 1. Therefore a tower $\sF$
 over $\F_q$ is 
asymptotically optimal if and only if $ \lambda (\sF) = \ell -1$  (for $  q = \ell ^2$). 
The tower $\sF = (F_0, F_1, F_2, \ldots)$ is called a {\it Galois tower} if all extensions
 $F_i /F_0$ are Galois. \\

From here on, $q = \ell ^2$ is a square. We will construct an asymptotically optimal Galois
 tower $\sE = (E_0, E_1, E_2, \ldots)$ over $\F_q$ with the properties stated in Theorem 1.7.
 The starting point is the asymptotically optimal tower $\sF = (F_0 , F_1, F_2, \ldots)$ 
over $\F_q$ which was introduced in [G-S2], see also [G-S3]. It is defined as follows: 
\begin{enumerate}
\item[i)] $F_0 = \F_q (x_0)$ is the rational function field.
 
\item[ii)] For all $i \geq 0$ we have $F_{i+1} = F_i (x_{i+1})$ with 
\end{enumerate}
$$
x^{\ell}_{i+1} + x_{i+1} = \frac{x^{\ell}_{i}}{x^{\ell -1}_{i} +1} . \eqno{(2.1)}
$$
We will need the following properties (F1) - (F5) of this tower $\sF$; see [G-S2, Sec. 3] for the
 proof of (F1), (F2), (F3), (F5), and [G-S3, Sec. 3] for the proof of (F4). 
\begin{enumerate}
\item[(F1)] All extensions $F_{i+1} /F_i$ are Galois of degree $\ell$. 
\item[(F2)] The only places of $F_0 = \F_q (x_0)$ which are ramified in the tower $\sF$, are
 the places $(x_0 = \alpha)$ with $\alpha^{\ell} + \alpha =0$ and the place $(x_0 = \infty)$. 
\item[(F3)] The places $(x_0 = \infty)$ and $(x_0 = \alpha)$ with $\alpha^{\ell -1} +1 =0$ are
 totally ramified in all extensions $F_n / F_0$; i.e., their ramification index in $F_n/F_0$ is $\ell ^n$. 
\item[(F4)] One can refine the extensions $F_{i+1} /F_i$ to Galois steps of degree 
$p = {\rm char} (\F_q)$ as follows: 
$$
F_i = H^{(0)}_{i} \subseteq H^{(1)}_{i} \subseteq \ldots \subseteq H^{(a)}_{i} = F_{i+1}
$$
with $[H^{(j+1)}_{i} : H^{(j)}_{i}] = p$. For any place $P$ of $H^{(j)}_{i}$ and $Q$ of
 $H^{(j+1)}_{i}$ lying above $P$, the different exponent $d(Q|P)$ satisfies 
$$
d(Q|P) = 2 (e(Q|P) -1). 
$$
\item[(F5)] All places $(x_0 = \alpha)$ of $F_0$ with $\alpha \in \F_q$ and 
$\alpha^{\ell}  + \alpha \neq 0$ split completely in the tower $\sF$; i.e., any of 
these places has $\ell ^n$ extensions in $F_n|F_0$, and all of them are rational 
places of $F_n$.
\end{enumerate}
\noindent  We set 
$$
w : = x^{\ell }_{0} + x_0 \ \ \mbox{and} \ \ z : = w^{\ell -1} ; \eqno{(2.2)}
$$
then 
$$
\F_q (z) \subseteq \F_q (w) \subseteq F_0 = \F_q (x_0) \subseteq F_1 \subseteq F_2 \subseteq \ldots \ . 
$$
The extension $\F_q (w) / \F_q (z)$ is cyclic of degree $(\ell -1)$, and the extension 
$F_0/ \F_q (w)$ is Galois of degree $\ell$. In the extension $F_0 / \F_q (z)$ we have 
the following ramification and splitting behaviour (which is easily checked): 
\begin{enumerate}

\item[(F6)] The place $(z = \infty)$ of $\F_q (z)$ is totally ramified in $F_0/\F_q (z)$; 
the only place of $F_0$ lying above $(z = \infty)$ is the place $(x_0 = \infty)$. 
\item[(F7)] Exactly $\ell$ places of $F_0$ lie above the place $(z =0)$, namely the
 places $(x_0 = \alpha)$ with $\alpha^{\ell} + \alpha =0$. Their ramification index in $F_0/ \F_q (z)$ is $\ell -1$. 
\item[(F8)] No other places of $\F_q (z)$ are ramified in $F_0$. 
\item[(F9)] One can refine the extension $F_0 / \F_q (w)$ to Galois steps of
 degree $p = {\rm char} (\F_q)$ as follows: 
$$
\F_q (w) = H^{(0)} \subseteq H^{(1)} \subseteq \ldots \subseteq H^{(a)} = F_0
$$
with $[H^{(j+1)} : H^{(j)}] = p$. For any place $P$ of $H^{(j)}$ and $Q$ of $H^{(j+1)}$
 lying above $P$, the different exponent $d(Q|P)$ satisfies 
$$
d(Q|P) =2 (e(Q|P) -1). 
$$
\item[(F10)] The place $(z=1)$ splits completely in the extension $F_0 /\F_q (z)$; the
 places of $F_0$ lying above $(z =1)$ are exactly the places $(x_0 = \alpha)$ with
 $\alpha \in \F_q$ and $\alpha^{\ell }+ \alpha \neq 0$. 
\end{enumerate}

After these preparations we can now prove Theorem 1.7. We start with the tower
 $\sF = (F_0, F_1, F_2, \ldots)$ as above; in particular we consider the elements 
$w, z \in F_0$ as defined in (2.2) above. Then we define the tower $\sE = (E_0, E_1, E_2,
 \ldots )$ as follows: $E_0 = \F_q (z)$ is the rational function field. For all $n \geq 1$, 
$$
E_n \ \mbox{is the Galois closure of field extension} \ F_{n-1}/E_0.
$$
We have then 
$$
E_0 = \F_q (z) \subseteq \F_q (w) \subseteq \F_q (x_0) \subseteq E_1 \subseteq E_2 \subseteq \ldots , 
$$
and items b), c) of Theorem 1.7 are clear. By Galois theory, the field $E_n$ is the composite of
 the fields 
$$
F_{n-1}, \tau (F_{n-1}), \rho (F_{n-1}), \ldots , 
$$
where $\tau, \rho, \ldots$ run through all embeddings of the field $F_{n-1}$ over $E_0$
 into a fixed algebraically closed field $\bar{E} \supseteq E_0$. The extension $\F_q (w)/E_0$
 is Galois, hence the field $\F_q (w)$ is mapped onto itself by all such embeddings
 of $F_{n-1}/E_0$. By items (F4) and (F9) above we can therefore obtain the field 
$E_n$ by iterated composites of $F_{n-1}$ with Galois extensions of degree
 $p = {\rm char} (\F_q)$. It follows that the degree of $E_n/F_q (w)$ is a power
 of $p$. Since $[F_{n-1} : \F_q (w)] = \ell ^n$, item d) of Theorem 1.7 follows.\\

We consider now the place $(z=1)$ of the rational function field $E_0 = \F_q (z)$.
 By items (F5) and (F10), this place splits completely in the extension $F_{n-1}/E_0$;
 hence it splits completely also in $\tau (F_{n-1}) /E_0$ for all embeddings $\tau$ as 
above. As follows from ramification theory, the place $(z=1)$ then splits completely in
 the composite field of $F_{n-1} , \tau (F_{n-1}), \ldots$ (see [S1, III.8.4]). We have thus
 proved item e) of Theorem 1.7. An immediate consequence is that $\F_q$ is the full constant
 field of $E_n$; this is item a) of Theorem 1.7. 
Item f) of Theorem 1.7 follows easily from (F3) and (F6).\\

 The core of the proof of Theorem 1.7 is item g). For 
its proof we need a result from [G-S3]: 

\begin{lemma}
Let $F/\F_q$ be a function field and let $G_1/F$ and $G_2 /F$ be linear disjoint Galois 
extensions of $F$, both of degree $p = {\rm char} (\F_q)$. Denote by $G = G_1 \cdot G_2$ 
the composite field of $G_1$ and $G_2$. Let $Q$ be a place of $G$ and denote by
 $Q_1, Q_2$ and $P$ its restrictions to the subfields $G_1, G_2$ and $F$. Suppose that we have 
$$
d(Q_i|P) = 2 (e(Q_i|P) -1), \ \mbox{for} \ i = 1,2. 
$$
Then $d(Q|Q_i) = 2 (e(Q|Q_i)-1)$ holds for $i = 1,2$. 
\end{lemma}

\begin{proof} 
See [G-S3, Lemma 1]. 
\end{proof}

Now we prove item f) of Theorem 1.7. First of all, it follows from items (F2), (F6), (F7), (F8) 
that the places $(w =0)$ and $(w=\infty)$ of $\F_q (w)$ are the only ramified places in
 $F_{n-1} /\F_q (w)$ and hence in $E_n/\F_q (w)$. We consider now a place $\tilde{Q}$
 of $E_n$ which is ramified in the extension $E_n/E_0$. By items (F2), (F6), (F7), 
(F8), $\tilde{Q}$ is either a zero or a pole of the function $w$; i.e., $\tilde{Q}$ is in 
the support of the divisor $A^{(n)}$ or $B^{(n)}$ (notation as in item f) of Theorem 1.7). \\

Let $Q : = \tilde{Q} \cap F_{n-1}$ be the restriction of $\tilde{Q}$ to the field $F_{n-1}$. We refine
 the extension $F_{n-1}/\F_q (w)$ to Galois steps of degree $p$: 
$$
\F_q (w) = K_0 \subseteq K_1 \subseteq \ldots \subseteq K_m = F_{n-1} \subseteq E_n, \eqno{(2.3)}
$$
with $[K_{j+1} :K_j] =p$. Let $P_j : = Q \cap K_j$ for $j = 0, \ldots , m$. By items (F4), (F9), the 
different exponents $d(P_{j+1} |P_j)$ are given by 
$$
d(P_{j+1} |P_j) =2 (e(P_{j+1} |P_j) -1), \ \mbox{for} \ j = 0, \ldots ,m-1. \eqno{(2.4)}
$$
The Galois closure $E_n$ of $F_{n-1}/E_0$ is obtained by iterated composites of the chain 
$$
K_0 \subseteq K_1 \subseteq \ldots \subseteq K_m
$$
with the chains 
$$
\tau (K_0 ) \subseteq \tau (K_1) \subseteq \ldots \subseteq \tau (K_m), 
$$
where $\tau$ runs through the embeddings of $F_{n-1} /E_0$. So we can refine the chain 
in (2.3) to a chain 
$$
\F_q (w) = K_0 \subseteq K_1 \subseteq \ldots \subseteq K_m = F_{n-1} \subseteq 
K_{m+1} \subseteq \ldots \subseteq K_r = E_n, 
$$
where all extensions $K_{j+1}/K_j$ are Galois of degree $p$ (for $ j =0, \ldots, r-1)$. 
We set $P_j : = \tilde{Q} \cap K_j$ for $j = m+1, \ldots, r$. It then follows from (2.4) 
and Lemma 2.1 that the different exponents $d(P_{j+1} |P_j)$ satisfy
$d(P_{j+1} |P_j) =2 (e(P_{j+1} |P_j)-1)$, for $j=0 , \ldots , r-1$. Using the transitivity of different
 exponents (cf. [S1, III.4.11]) we obtain that 
$$
d (\tilde{Q} |P_0) =2 (e(\tilde{Q} |P_0) -1). 
$$
This finishes the proof of item g) of Theorem 1.7. \\

With notations as in items f) and g), the Hurwitz genus formula for the extension
 $E_n /\F_q (w)$ yields
$$
\begin{array}{ll} 
2g (E_n) -2 & = -2 [E_n : \F_q (w)] + 2e^{(n)}_{0} \deg A^{(n)} \\
& \\
 &+ 2 e^{(n)}_{\infty} \deg B^{(n)} -2 (\deg A^{(n)} + \deg B^{(n)} ) \\
& \\
& = 2 \cdot [E_n : \F_q (w)] -2 (\deg A^{(n)} + \deg B^{(n)} ). 
\end{array}
$$
We have used here that the divisors $e^{(n)}_{\infty} \cdot B^{(n)}$ and $e^{(n)}_{0} \cdot A^{(n)}$ 
are the pole divisor and the zero divisor of the function $w$ in $E_n$, hence their degree is
 equal to the degree $[E_n : \F_q (w)]$. We have thus proved item h) of Theorem 1.7. \\

From items e) and h) we see that 
$$
N(E_n) /g(E_n) \geq \ell -1 \ \mbox{for all} \ n \geq 1, \eqno{(2.5)}
$$
Hence $\lim_{n \to \infty} N(E_n)/g(E_n) \geq \ell -1$. By the Drinfeld-Vladut bound (see Section  1) we
 also have that $\lim_{n \to \infty} N(E_n)/g(E_n) \leq \ell -1$, hence equality holds. This proves
 item j) and finishes the proof of Theorem 1.7. 

\section{Asymptotically Good Transitive Codes}

The aim of this section is to prove Theorem 1.5. We use notation as in Theorem 1.5, in
 particular $q = \ell ^2$ is a square. Let $R, \delta \geq 0$ with 
$$
R = 1 - \delta - \frac{1}{\ell -1} , \eqno{(3.1)}
$$
and let $\epsilon > 0$. We will construct transitive codes $C$ over $\F_q$ of arbitrarily large
 length such that $R(C) \geq R - \epsilon$ and $\delta (C) \geq \delta$; this proves then Theorem 1.5. \\

Consider the tower $\sE = (E_0, E_1, E_2, \ldots )$ of function fields over $\F_q$ which 
was constructed in Theorem 1.7. Choose an integer $n > 0$ so large that 
$$
\frac{1}{\ell ^n (\ell -1)} < \epsilon. \eqno{(3.2)}
$$
Let $N: = [E_n : \F_q (z)]$, with the function $z \in E_n$ as in Theorem 1.7, and consider 
the divisors $D, G_0$ of $E_n$ which are given by 
$$\\
D: = \sum_{P|(z=1)} P \ \mbox{and} \ G_0 : = \sum_{Q|(z=\infty)} Q. \eqno{(3.3)}
$$
This means: $P$ runs over all places of $E_n$ which are zeroes of the function $(z-1)$,
and $Q$ runs over all poles of the function $z$ in $E_n$.
By Theorem 1.7 e) all these places $P$ are rational, and the degree of $D$ is $\deg D =N$.
 With notations as in Theorem 1.7 f), the divisor $G_0$ is just the divisor $G_0 = B^{(n)}$,
 since the functions $w$ and $z = w^{\ell -1}$ have the same poles. The degree of $G_0$
 satisfies then 
$$
\deg G_0 = \frac{[E_n : \F_q (w)]}{e^{(n)}_{\infty}} \leq \frac{[E_n : \F_q (w)]}{\ell ^n} = \frac{N}{\ell ^n (\ell -1)} , 
$$
by Theorem 1.7 f). Hence we have that 
$$
(\deg G_0) / N < \epsilon, 
$$
by Inequality (3.2). We choose  $r \geq 0$ such that 
$$
1 - \delta \geq r \cdot \frac{\deg G_0}{N} > 1 - \delta - \epsilon \eqno{(3.4)}
$$
and consider the {\it geometric Goppa code} 
$$
C: = C_{\sL} (D, r G_0) \subseteq \F^{N}_{q}
$$
associated to the divisors $D$ and $rG_0$. It is defined as follows (cf. [S1, II.2.1] or [T-V]): 
If $\sL (rG_0) \subseteq E_n$ denotes the Riemann-Roch space of the divisor $rG_0$ 
and the divisor $D$ is defined as $D = P_1 + \ldots +P_N$, then 
$$
C_{\sL} (D,rG_0) = \{ (f(P_1), \ldots , f(P_N)) \in \F^{N}_{q} | f \in \sL (rG_0) \} . \eqno{(3.5)} 
$$
For the parameters $k = \dim C$ and $d =d (C)$ we have the standard estimates for 
geometric Goppa codes (see [S1, II.2.3]): 
$$
k \geq r \cdot \deg G_0 +1 - g (E_n) \ \ \mbox{and} \ \ d \geq N - r \cdot \deg G_0.
$$
Hence the information rate $R (C)$ satisfies 
$$
R(C) = \frac{k}{N} \geq \frac{r \cdot \deg G_0}{N} + \frac{1}{N} - \frac{g(E_n)}{N} >
 1 - \delta - \epsilon - \frac{g(E_n)}{N} \ ,
$$
by Inequality (3.4). Now observe that 
$$
\frac{g(E_n)}{N} \leq \frac{1}{\ell -1}, 
$$
by Inequality (2.5), and we obtain using Equality (3.1) the following estimate for $R(C)$: 
$$
R(C) > 1 - \delta - \epsilon - \frac{1}{\ell -1} = R - \epsilon. 
$$
For the relative minimum distance $\delta (C)$ we get with (3.4): 
$$
\delta (C) = \frac{d}{N} \geq \frac{N - r \cdot \deg G_0}{N} = 1 - \frac{r \cdot \deg G_0}{N} \geq \delta. 
$$
These are the desired inequalities for $R(C)$ and $\delta (C)$.\\

 It remains to show 
that the code $C= C_{\sL} (D, r G_0)$ that we constructed above is in fact a {\it transitive}
 code. To this end we consider the Galois group of the extension $E_n/E_0$, 
$$
\Gamma : = {\rm Gal} (E_n /E_0).
$$
The places $P_1, \ldots , P_N$ in the support of the divisor $D$ are exactly the
 places of $E_n$ lying above the place $(z =1)$; hence $\Gamma$ acts transitively 
on the set $\{ P_1, \ldots , P_N\}$, see [S1, III.7.1]. The divisor $r G_0$ is obviously
 invariant under the action of $\Gamma$. Therefore $\Gamma$ acts on the code
 $C = C_{\sL} (D, r G_0)$ as a transitive permutation group in the following way
 (see [S1, VII.3.3]): for $\sigma \in G$ and $f \in \sL (r G_0)$, 
$$
\sigma (f(P_1), \ldots , f(P_N)) = (f (\sigma P_1), \ldots , f(\sigma P_N)). 
$$
This completes the proof of Theorem 1.5 \qed 

\begin{remark}
It is an obvious idea to prove the existence of asymptotically good {\it cyclic} 
codes in a similar manner. One should start with a tower $\sH = (H_0, H_1,
 H_2, \ldots)$ of function fields over $\F_q$, where all extensions $H_n/H_0$ 
are cyclic Galois extensions; then one can do the same construction of codes 
as in the proof of Theorem 1.5 above. However, this method does not work: it is
 known that the limit $\lambda (\sH) = \lim_{n \to \infty} N(H_n)/g(H_n)$ of such 
a ``cyclic'' tower $\sH$ is zero, see [F-P-S]. 
\end{remark}

\begin{remark}
The notion of ``information rate''of a code can be defined also for {\it non-linear codes} 
$C \subseteq \F^{N}_{q}$, by setting $R(C): = \log_q (|C|/N)$. Using this definition,
 one obtains in an obvious manner an analogue of the function $\alpha_q (\delta)$ 
by considering {\it all} codes over $\F_q$, not just linear codes. We denote this 
analoguous function again by $\alpha_q (\delta)$. It was shown in [N-O1] and [S-X] 
that in a large open subinterval of [0,1], the Tsfasman-Vladut-Zink bound
$$
\alpha_q (\delta) \geq 1 - \delta - A(q)^{-1} \eqno{(3.6)}
$$
can be improved to 
$$
\alpha_q (\delta) \geq 1 - \delta - A (q)^{-1} + \log_q (1 + q^{-3}). \eqno{(3.7)}
$$
\end{remark}

A further slight improvement of Inequality (3.7) was very recently found in [N-O2]. However, 
it seems that the codes which were constructed in [N-O1,2] and [S-X] in order to prove 
Inequality (3.6) do not have any algebraic or combinatoric structure. By combining the 
method of [S-X] with our proof of Theorem 1.5 we can now show that the lower bound (3.7) 
for $\alpha_q (\delta)$ is attained by {\it transitive} non-linear codes. 

\begin{theorem}
Assume that $q = \ell ^2$ is a square, and set 
$$\delta^* : = 1 - 2 /(\ell -1) - (4q -2)/((q-1) (q^3 +1)).$$
 Then the bound
$$
\alpha_q (\delta) \geq 1 - \delta - A(q)^{-1} + \log_q (1 + q^{-3})
$$
is attained by transitive codes, for all $\delta$ in the interval $(0, \delta^*) \subseteq [0,1]$. 
\end{theorem}

\begin{proof} 
(Sketch). We recall briefly the code construction given in [S-X]. One considers a
 function field $F$ over $\F_q$ of genus $g$ and a set $\sP = \{ P_1, \ldots , P_N\}$ 
of $N$ distinct rational places of $F$. Let $H \geq 0$ be a divisor of $F$ of degree 
$\deg H \geq 2g -1$ with ${\rm supp} (H) \cap \sP = \emptyset$ and consider divisors
 $G$ of the form 
$$
G = \sum^{t}_{j=1} m_{i_j} P_{i_j} \ \mbox{with} \ 1 \leq i_1 < i_2 < \ldots < i_t \leq N, 
m_{i_j} \geq 1 \ \mbox{and} \ \deg G = s. \eqno{(3.8)}
$$
Define the set $M_H (G)$ as follows: 
$$
M_H (G) : = \{ x \in \sL (H+G)  \ | \ v_{P_{i_j}} (x) = - m_{i_j} \ \mbox{for} \ 1 \leq j \leq t \} .
$$
Choose integers $s,t$ with $1 \leq t \leq N$ and $s \geq t$, and set 
$$
S : = S (H, \sP, s, t) : = \bigcup_G M_H (G), 
$$
where $G$ runs over all divisors of the form (3.8). It is clear that $M_H (G_1) \cap
 M_H (G_2) = \emptyset$ if $G_1 \neq G_2$. Hence we can define a map $\varphi :
 S \to \F^{N}_{q}$ in the following way: for $x \in M_H (G)$ put $\varphi (x) = (x_1, \ldots , x_N)$ with 
$$
x_i = \left\{ \begin{array}{cl} x(P_i) & \mbox{if} \ P_i \not\in {\rm supp} (G), \\
0 & \mbox{if} \ P_i \in {\rm supp} (G). \end{array} \right. 
$$
Thus we obtain a (non-linear) code $C= C (H, \sP, s,t)$ by setting 
$$
C(H, \sP, s, t): = \varphi (S) \subseteq \F^{N}_{q}.
$$
If the function field $F$ runs through a sequence of function fields $(F_0, F_1, F_2, 
\ldots)$ over $\F_q$ with $\lim_{n \to \infty} N(F_n)/g(F_n) = \sqrt{q} -1$, one can 
choose the set $\sP$, the divisor $H$ and the integers $s,t$ in such a way that 
the corresponding codes $C(H, \sP, s, t)$ reach the bound (3.7), see [S-X, Prop.3.3 
and Thm.3.4.] \\
 
In order to obtain transitive codes with the above construction, we use again the function
 fields $E_n$ of the tower $\sE = (E_0, E_1, E_2, \ldots)$ from Theorem 1.7. We choose the 
set $\sP$ as in the proof of Theorem 1.5; i.e., 
$$
\sP = \{ P \ | \  P \ \mbox{is a zero of the function} \ z -1 \ \mbox{in} \ E_n \} , 
$$
see (3.3). The divisor $H$ is chosen as 
$$
H = m_0 \cdot G_0,
$$
with the divisor $G_0$ of $E_n$ as in (3.3). Since the set $\sP$ and the divisor $G_0$ 
are invariant under the action of the group $\Gamma = {\rm Gal} (E_n/E_0)$, it follows 
immediately that the corresponding codes $C(H,\sP, s,t) \subseteq \F^{N}_{q}$ are 
$\Gamma$-invariant; they are therefore transitive codes. 
\end{proof} 

\section{Asymptotically Good Self-Dual and Self-Orthogonal Codes}

In this section we shall prove Theorem 1.6. First we recall some definitions
and facts. 

\begin{definition}
Let $C \subset \F^{N}_{q}$ be a linear code, and let $\underline{a} = (a_1, \ldots , a_N) \in \F^N_q$
 with non-zero components $a_1, \ldots , a_N \neq 0$. We set 
$$
\underline{a} \cdot C : = \{ (a_1 \cdot c_1, \ldots , a_N \cdot c_N) \in \F^{N}_{q} \ | \ (c_1, \ldots , c_N) \in C \} , 
$$
and call the codes $C$ and $\underline{a} \cdot C$ {\it equivalent}. 
\end{definition}

It is clear that equivalent codes have the same parameters (length, dimension, minimum 
distance, information rate, relative minimum distance). Note however that the automorphism 
groups ${\rm Aut} (C)$ and ${\rm Aut} (\underline{a} \cdot C)$ are in general non-isomorphic. 

\begin{definition}
i) A code $C \subseteq \F^{N}_{q}$ is called {\it self-dual} if $C$ is equal to its dual code
 $C^{\perp}$. The code $C$ is called {\it self-orthogonal} if $C \subseteq C^{\perp}$.\\
ii) A code $C$ is called {\it iso-dual} if $C$ is equivalent to its dual code $C^{\perp}$, cf. [P-H]. \\
iii) A code $C$ is called {\it iso-orthogonal} if $C$ is equivalent to a subcode of $C^{\perp}$. 
\end{definition}

Now let $F/\F_q$ be a function field and let $P_1 , \ldots , P_N$ be distinct rational places
 of $F$. Let $D = P_1 + \ldots + P_N$ and let $G$ be a divisor with ${\rm supp} \ D
 \cap {\rm supp} \ G = \emptyset$. As in Section 3, we consider the {\it geometric Goppa code} (cf. (3.5))
$$
C_{\sL} (D, G) : = \{ (f(P_1), \ldots , f(P_N)) \in \F^{N}_{q} \ | \ f \in \sL (G) \}. \eqno{(4.1)}
$$

\begin{proposition}
Let $D$ and $G$ be divisors of the function field $F/\F_q$ as above and consider
 the code $C = C_{\sL} (D, G)$ as defined in {\rm (4.1)}. Suppose that $\eta$ is a differential
 of $F$ with the property $v_{P_i} (\eta ) = -1$ for $i =1, \ldots , N$. Then the dual code
 $C^{\perp} = C_{\sL} (D, G)^{\perp}$ is given by 
$$
C^{\perp} = \underline{a} \cdot C_{\sL} (D,H), 
$$
with $H: = D - G + (\eta)$ and $\underline{a} = ({\rm res}_{P_1} (\eta) , \ldots , {\rm res}_{P_N} (\eta ))$. 
\end{proposition}

\begin{proof}
See [S1, Cor.2.7]. 
\end{proof}

We want to apply Proposition 4.3 to geometric Goppa codes which are defined by means of the
 function fields $E_n$ in the tower $\sE = (E_0, E_1, E_2, \ldots)$ of Theorem 1.7. So we
 must find an appropriate differential $\eta$ of $E_n$ having the properties as required in Proposition 4.3.

\begin{proposition}
We assume all notations from Theorem {\rm 1.7} and consider the differential 
$$
\eta : = \frac{dw}{1-z} 
$$
of the function field $E_n$ (with $n \geq 2$). Then the following holds: 
\begin{enumerate}
\item[i)] The divisor of $\eta$ in $E_n$ is given by 
$$
(\eta) = a_n \cdot A^{(n)} + b_n \cdot B^{(n)} - D^{(n)}, 
$$
where the divisors $A^{(n)} > 0$ and $B^{(n)} > 0$ are as in Theorem {\rm 1.7 f) }, the integers 
$a_n > 0$ and $b_n > 0$ satisfy $a_n \equiv b_n \equiv 0 \ {\rm mod} \ 2$, and the
 divisor $D^{(n)}$ is the sum over all zeroes of the function $z -1$ in $E_n$; i.e., 
$$
D^{(n)} = \sum_{P | (z=1)} P. 
$$
\item[ii)] The residue of the differential $\eta$ at a place $P$, which is a zero of $z-1$ in 
$E_n$, is an element of $\F^{\times}_{\ell}$; i.e., 
$$
{\rm res}_P (\eta ) = \alpha_P \ \ \mbox{with} \ \ \alpha^{\ell -1}_{P} =1.
$$
\end{enumerate}
\end{proposition}

\begin{proof}
i) By Theorem 1.7 f), the principal divisor of the function $w$ in $E_n$ is 
$$
(w)^{E_n} = e^{(n)}_{0} \cdot A^{(n)} - e^{(n)}_{\infty} \cdot B^{(n)}, 
$$
and by item g) of Theorem 1.7, the different of $E_n /\F_q (w)$ is 
$$
{\rm Diff} (E_n / \F_q (w)) = 2 (e^{(n)}_{0} -1) \cdot A^{(n)} + 2 (e^{(n)}_{\infty} -1) \cdot B^{(n)}. 
$$
It follows that the divisor of the differential $dw$ in $E_n$ is given by (see [S1, III.4.6]) 
$$
(dw)  = -2 e^{(n)}_{0} B^{(n)} + {\rm Diff} (E_n /\F_q (w)) 
 = 2e^{(n)}_{0} A^{(n)} - 2 A^{(n)} - 2 B^{(n)} .
$$
The divisor of the function $1 - z$ in $E_n$ is 
$$
(1 -z)^{E_n} = D^{(n)} - (\ell -1) \cdot e^{(n)}_{\infty} \cdot B^{(n)}, 
$$
and we obtain the divisor of the differential $\eta = {dw}/{(1-z)}$ as follows: 
$$ 
\begin{array}{ll}
(\eta) & = 2e^{(n)}_{0} A^{(n)} - 2 A^{(n)} - 2 B^{(n)} - D^{(n)} + (\ell -1) e^{(n)}_{\infty} B^{(n)} \\
& = a_n A^{(n)} + b_n B^{(n)} - D^{(n)} , 
\end{array} 
$$
with $a_n>0, \ b_n>0$ and $a_n \equiv b_n \equiv 0 \mod 2$. \\

\noindent ii) Let $P$ be a place of $E_n$ which is a zero of the function $z -1$. The element $t: = z -1$ is
 a $P$-prime element. From the equation $w^{\ell -1} = z = t +1$ we obtain 
$$
dt = (\ell -1) w^{\ell -2} dw = - \frac{w^{\ell -1}}{w} dw = - \frac{1+t}{w} dw, 
$$
hence 
$$
\eta = \frac{dw}{1-z} = - \frac{1}{t} dw = \frac{w}{1+t} \cdot \frac{1}{t} dt. 
$$
Let $\alpha : = w (P) \in \F_q$ be the residue class of $w$ at the place $P$; then 
$$
\frac{w}{1+t} \equiv \alpha \ {\rm mod} \ P \ \mbox{and therefore} \ {\rm res}_P (\eta) = \alpha. 
$$
Since $\alpha^{\ell -1} = w^{\ell -1} (P) = z (P) =1$, we conclude that $\alpha \in \F_{\ell} \setminus \{ 0 \}$. 
\end{proof} 

Now we can construct certain geometric Goppa codes which are associated to the function
 field $E_n$ in the tower $\sE = (E_0, E_1, E_2, \ldots)$ of Theorem 1.7. For the rest of this 
section we fix notations as above; in particular we will use without further explanation the
 divisors $A^{(n)}, B^{(n)}$ and $D^{(n)}$, the differential $\eta$ and the integers $a_n$ 
and $b_n$ as in Proposition 4.4. 

\begin{definition}
For integers $a,b$ with $0 \leq a \leq a_n$ and $0 \leq b \leq b_n$, we define the 
code $C^{(n)}_{a,b}$ by 
$$
C^{(n)}_{a,b} : = C_{\sL} (D^{(n)} , a A^{(n)} + b B^{(n)}). 
$$
\end{definition}

\begin{remarks} i) It is clear that the codes $C^{(n)}_{a,b}$ are transitive. This follows
 as in Section 3 from the fact that the Galois group $\Gamma = {\rm Gal} (E_n /E_0)$ acts 
transitively on the places $P \in {\rm supp} (D^{(n)})$ and leaves the divisors $A^{(n)}$
 and $B^{(n)}$ invariant. \\[.1cm]
ii) For $n \to \infty$, the codes $C^{(n)}_{a,b}$ attain the Tsfasman-Vladut-Zink 
bound $\alpha_q (\delta) \geq 1 - \delta -1 /(\ell -1)$, for all $\delta \in (0, 1 -1/(\ell -1))$. 
This is proved in the same manner as Theorem 1.5 (see Section 3). 
\end{remarks}

\begin{proposition}
We write $D^{(n)} = P_1 + \ldots + P_N$, with $N = [E_n : E_0]$, and set 
$$
\underline{u} : = ({\rm res}_{P_1} \eta , \ldots , {\rm res}_{P_N} \eta ) \in (\F^{\times}_{q} )^N.
$$
Then the dual of the code $C^{(n)}_{a,b}$ is given by 
$$
(C^{(n)}_{a,b} )^{\perp} = \underline{u} \cdot C^{(n)}_{a_n -a, b_n -b} . 
$$
\end{proposition}

\begin{proof} 
The differential $\eta$ satisfies the condition $v_{P_i} (\eta) =-1$, for $i =1, \ldots ,N$. Hence
 we can apply Proposition 4.3 and obtain 
$$
(C^{(n)}_{a,b})^{\perp} = \underline{u} \cdot C_{\sL} (D^{(n)} , H), 
$$
with 
$$
\begin{array}{ll}
H & = D^{(n)} - (a A^{(n)} + b B^{(n)} ) + (\eta) \\
& \\
& = D^{(n)} - (a A^{(n)} + b B^{(n)}) + (a_n A^{(n)} + b_n B^{(n)} - D^{(n)} ) \\
& \\
& = (a_n -a) A^{(n)} + (b_n -b) B^{(n)} . 
\end{array}
$$
We have used here Proposition 4.4 i). 
\end{proof}

The following corollary is an obvious consequence from Proposition 4.7, cf. Definition 4.2. 

\begin{corollary}
{\rm i)} For $0 \leq a \leq a_n /2$ and $0 \leq b \leq b_n /2$, the code $C^{(n)}_{a,b}$ is 
transitive and iso-orthogonal. \\
{\rm ii)} For $a = a_n /2$ and $b = b_n /2$, the code $C^{(n)}_{a,b}$ is iso-dual. 
\end{corollary}

\begin{corollary}
{\rm i)} For $0 \leq a \leq a_n/2$ and $0 \leq b \leq b_n /2$, the code $C^{(n)}_{a,b}$ is 
equivalent to a self-orthogonal code $\tilde{C}^{(n)}_{a,b}$. \\
{\rm ii)} For $a = a_n /2$ and $b= b_n /2$, the code $C^{(n)}_{a,b}$ is equivalent to a 
self-dual code $\tilde{C}^{(n)}_{a,b}$. 
\end{corollary}

\begin{proof}
The components of the vector $\underline{u} = ({\rm res}_{P_1} \eta, \ldots , {\rm res}_{P_N} \eta)$ 
in Proposition 4.7
 are in $\F^{\times}_{\ell }$, by Proposition 4.4 ii). So we can write ${\rm res}_{P_i} \eta = v^{2}_{i}$ with
 $v_i \in \F^{\times}_{q}$ (note that $q = \ell ^2$). We set $\underline{v} : = (v_1, \ldots , v_N)$; then the code 
$$
\tilde{C}^{(n)}_{a,b} : = \underline{v} \cdot C^{(n)}_{a,b}
$$
is self-orthogonal, resp. self-dual. 
\end{proof} 

Theorem 1.6 is now an immediate consequence of Corollary 4.9 and Remark 4.6 ii). 

\begin{remark}
The existence of asymptotically good sequences $(C_j)_{j \geq 0}$ of isodual geometric 
Goppa codes over $\F_q$ (with $q = \ell ^2$) was already proved in [Sch]. However, the codes 
which were constructed there attain only the lower bound 
$$
\lim_{j \to \infty} \delta (C_j) \geq \frac{1}{2} - \frac{1}{\ell -3} . \eqno{(4.2)}
$$
The codes $\tilde{C}_n : = \tilde{C}^{(n)}_{a,b}$ in Corollary 4.9 ii) are not only iso-dual but they 
are self-dual. They satisfy the bound (see Theorem 1.6 ii))
$$
\lim_{j \to \infty} \delta (\tilde{C}_j) \geq \frac{1}{2} - \frac{1}{\ell -1}, 
$$
which is better than Inequality (4.2). 
\end{remark}

\section{Conclusion}

Let $q = \ell ^2$ be a square. We have shown in this paper, that the following classes
 of linear codes over $\F_q$ attain the Tsfasman-Vladut-Zink bound: 
\begin{enumerate}
\item[-] {\it transitive} codes (Theorem 1.5),\\
 
\item[-] {\it transitive iso-orthogonal} codes (Corollary 4.8),\\
 
\item[-] {\it transitive iso-dual} codes (Corollary 4.8),\\
 
\item[-] {\it self-orthogonal} codes (Theorem 1.6), \\

\item[-] {\it self-dual} codes (Theorem 1.6).\\
 
\end{enumerate}

In particular, the above classes of codes are better than the Gilbert-Varshamov 
bound, for all squares $q \geq 49$.  The class of {\it non-linear transitive} codes 
attains an even better bound than the Tsfasman-Vladut-Zink bound (Theorem 3.3).

\bibliographystyle{plain}

\end{document}